\documentclass[10pt]{article}

\usepackage{amsthm}
\usepackage{amssymb}
\usepackage{amsmath}

\usepackage[all]{xy}
\usepackage{graphicx}

\newtheorem{thm}{Theorem}[section]
\newtheorem{cor}[thm]{Corollary}
\newtheorem{lem}[thm]{Lemma}
\newtheorem{prop}[thm]{Proposition}

\newtheorem{defn}[thm]{Definition}

\theoremstyle{definition}
\newtheorem{ex}{Example}[section]

\theoremstyle{remark}

\newcommand{\secref}[1]{\S\ref{#1}}

\newcommand{\mA}{\mathcal A}
\newcommand{\mB}{\mathcal B}
\newcommand{\mC}{\mathcal C}
\newcommand{\mE}{\mathcal E}

\newcommand{\mG}{\mathcal G}
\newcommand{\mH}{\mathcal H}

\newcommand{\mL}{\mathcal L}
\newcommand{\mM}{\mathcal M}
\newcommand{\mO}{\mathcal O}

\newcommand{\mZ}{\mathcal Z}

\newcommand{\mUE}{\mathcal {UE}}
\newcommand{\wE}{\wa \mE}

\newcommand{\coe}{\mO_{\mE}}

\newcommand{\com}{\mO_{\mM}}

\newcommand{\zro}{\mZ^\rho}
\newcommand{\oro}{\mO_{\rho}}
\newcommand{\orok}{\mO_\rho^k}
\newcommand{\soro}{\mO_{\left[ \rho \right]}}

\newcommand{\wa}{\widehat}

\newcommand{\ra}{\rightarrow}

\newcommand{\hra}{\hookrightarrow}

\newcommand{\sC}{{\it C*}-}

\newcommand{\eps}{\varepsilon}

\newcommand{\bC}{\mathbb C}

\newcommand{\bN}{\mathbb N}
\newcommand{\bS}{\mathbb S}
\newcommand{\bT}{\mathbb T}
\newcommand{\bZ}{\mathbb Z}
\newcommand{\bP}{\mathbb P}

\newcommand{\sud}{{\mathbb {SU}}(d)}

\newcommand{\ii}{\iota,\iota}
\newcommand{\rr}{\rho,\rho}
\newcommand{\rs}{\rho,\sigma}

\newcommand{\sss}{\sigma,\sigma}

\newcommand{\cpen}{\mA \rtimes^\mE_\rho \bN }
\newcommand{\cpn}{\mA \rtimes_\rho \bN }

\begin{document}

\author{Ezio Vasselli 
                         \\Dipartimento di Matematica
                         \\University of Rome "Tor Vergata" 
                         \\E-mail {\tt vasselli@mat.uniroma2.it} }

\title{Crossed products of $C(X)$-algebras by endomorphisms and $C(X)$-categories. 
       \\
       ({\it Preliminary Draft})}

\maketitle

\begin{abstract}
We construct the crossed product of a $C(X)$-algebra by an endomorphism, in such a way that it becomes induced by a Hilbert $C(X)$-bimodule. Furthermore we introduce the notion of $C(X)$-category, and discuss relationships with crossed products and duality for compact groups.
\end{abstract}

\section{Introduction.}

Basic works by S. Doplicher and J.E. Roberts (\cite{DR87,DR88,DR89,DR89A}) established at the end of the eighties a deep link between duality for compact groups and theory of crossed products. In that context, the main result is characterization of abstract strict tensor \sC categories (in the sense of \cite{DR89}, see also \secref{categories} in the present paper) that are isomorphic to the dual of a compact group. We roughly recall the main steps of the proof, without pretending to be exhaustive, and assuming for simplicity that our strict tensor \sC category $\wa \rho$ has objects the tensor powers $\rho^r$ of a unique object $\rho$. Thus we can label objects simply by natural numbers: $r \sim \rho^r , 0 \sim \iota$, where $\iota$ is the identity object. The symbol $1_\rho$ will denote the identity arrow.

\begin{itemize}

\item As a first step a unital \sC algebra (called DR-algebra in the following) $\oro$ is canonically associated to the object $\rho$, by identifying arrows $t \sim t \times 1_\rho$, $t \in (r,s)$ (here $\times$ denotes the tensor product) and assigning suitable multiplication (induced by composition of arrows), involution and norm on the quotient (see \cite{DR89}). If $(0,0) \simeq \bC$ then $\oro$ has trivial centre. Well known \sC algebras arise in this way, as Cuntz-Krieger-Pimsner algebras in the case in which $\rho$ is a Hilbert \sC bimodule; $(r,s)$ is then the Banach space of bimodule maps from the $r$-fold tensor power of the bimodule into the $s$-fold tensor power (see \cite{DPZ97} for details).

\item $\oro$ comes equipped of a canonical endomorphism $\sigma_\rho$, induced by left tensoring by the identity arrow: $\sigma_\rho (t) := 1_\rho \times t , t \in (r,s)$. Under suitable hypotesis, the identity $(r,s) = (\sigma_\rho^r , \sigma_\rho^s)$ holds (this property is called \sC amenability in \cite{LR97}). Thus $\wa \rho$ is characterized as the category of intertwiners in $\oro$ of powers of the canonical endomorphism $\sigma_\rho$. Properties of $\wa \rho$ as symmetry of tensor product and existence of conjugates can be translated as properties of the endomorphism $\sigma_\rho$. These properties are respectively permutation symmetry with order some $d \in \bN$, and special conjugate property (see \cite{DR87,DR89A}).

\item A basic crossed product construction can now be applied (\cite{DR89A}). Let $\tau$ be an endomorphism of a unital \sC algebra $\mA$ with trivial centre, carrying permutation symmetry of order $d$ and satisfying the special conjugate property. Then there is a unique \sC algebra $\mA \rtimes \wa \tau$ satisfying the following properties:

   \begin{itemize}
 
   \item There is a unital monomorphism $i : \mA \hra \mA \rtimes \wa \tau$, and a Hilbert space $\mH := {\bf span} \left\{ \psi_i \right\}_{i=1,..,d} \subset \mA \rtimes \wa \tau$ of orthogonal isometries with total support $1$ (see \cite{DR88} about the notion of Hilbert space in a \sC algebra). The \sC subalgebra generated by $\mH$ is thus isomorphic to the Cuntz algebra $\mO_d$. The relative commutant of $i(\mA)$ in $\mA \rtimes \wa \tau$ is trivial;
 
   \item $\mH$ induces canonically an endomorphism $\sigma_\mH (b) := \sum_i \psi_i b \psi_i^*$ on $\mA \rtimes \wa \tau$, and $i \circ \tau = \sigma_\mH \circ i$;
 
   \item a unique subgroup $G \subseteq \sud$ acts by automorphisms on $\mA \rtimes \wa \tau$. The fixed point algebra is $i(\mA)$. The $G$-action coincides with the natural action on $\mH \simeq \bC^d$ and extends in the usual way to the Cuntz algebra generated by $\mH$.

   \end{itemize} 
 
Now in this setting it is easy to prove that $(\tau^r,\tau^s)$ coincides with the intertwiners space $(\mH^r,\mH^s)_G$ of tensor powers of the defining representation of $G$. Thus in the case $\mA = \oro$, $\tau = \sigma_\rho$ we get by the previous points $(\rho^r,\rho^s) \simeq (\mH^r,\mH^s)_G$, and $\wa \rho$ is isomorphic to the tensor category of powers of the defining representation of $G$.
 
\end{itemize}

\bigskip
\bigskip

Aim of the present paper is to expose some preliminary results (that are part of PhD thesis of the author) as a starting point to extend the above sketched construction to the case in which the abelian \sC algebra $(\ii)$ does not reduce to the complex numbers. We will prove in the following sections that further structure in the more general case naturally appears, at crossed product and categorical level. This paper is organized as follows:

\bigskip

In \secref{hilb_bimod} we recall some terminology and basic facts about Hilbert bimodules in \sC algebras, and relationship with endomorphisms. Furthermore a generalized notion of inner endomorphism is introduced.

\bigskip

In \secref{cp_vec_bundle} we construct the crossed product of a $C(X)$-algebra by a $C(X)$-endomorphism. In the crossed product \sC algebra the endomorphism becomes induced by the module of continuous sections of a vector bundle over $X$. This setting is natural, in the sense that every endomorphism of a \sC algebra induces a structure of $C(X)$-algebra, and generalizes the crossed products by Cuntz \cite{Cun82}, Paschke \cite{Pas80}, Stacey \cite{Sta93}. An important tool (and example) is the Pimsner algebra associated to the module of continuous sections of the vector bundle.

\bigskip

In \secref{categories} we introduce the notion of $C(X)$-category, and give a characterization in terms of 'fields of \sC categories'. Furthermore we prove that every strict tensor \sC category has a natural structure of $C(X)$-category, and give a description of DR-algebras in terms of fields of \sC algebras. In particular, the strict tensor \sC category with arrows intertwiners of powers of a \sC algebra endomorphism is an example of $C(X)$-category.

\bigskip

In \secref{dual_dyn_sys} we introduce the notion of cross product by a dual action, that is the natural framework for duality applied to the strict tensor \sC category $\wa \tau$ of intertwiners of a \sC algebra endomorphism $\tau$. In this setting we give an algebraic condition to be $\wa \tau$ isomorphic to the strict tensor $C(X)$-category of tensor powers of a faithful compact group representation over a suitable Hilbert bimodule.

\section{Hilbert Bimodules in \sC algebras.}
\label{hilb_bimod}

We recall some definitions from \cite{DPZ97}. Let $\mB \subset \mA$ be an inclusion of \sC algebras. A {\em Banach} $\mB${\em -bimodule in} $\mA$ is a closed vector space $\mM \subset \mA$ such that $\mB \cdot \mM , \mM \cdot \mB \subseteq \mM$. In particular $\mM$ is said to be {\em Hilbert} $\mB${\em -bimodule in} $\mA$ if $\mM^* \cdot \mM \subseteq \mB$, in such a way that the map $\psi , \psi' \mapsto \psi^* \psi'$ defines a $\mB$-valued scalar product on $\mM$.

In the theory of Hilbert bimodules it is customary to require that left action of the coefficients \sC algebra is non degenerate. This fact can be encoded for Hilbert bimodules in \sC algebras by requiring that the left annihilator of $\mM$ in $\mB$ is trivial.

Banach and Hilbert bimodules appear naturally in \sC algebras by considering endomorphisms. Let in fact $\mZ$ denote the centre of $\mA$, $\rho$ an endomorphism on $\mA$. Then the spaces of intertwiners $(\rho^r , \rho^s) := \left\{ t \in \mA : \rho^s(a) t = t \rho^r(a) , a \in \mA \right\}$ are Banach $\mZ$-bimodules. In particular for $r = 0$ (as usual we define $\rho^0$ to be the identity automorphism $\iota$) we find that $(\iota , \rho^s)^* \cdot (\iota , \rho^s) \subseteq (\ii) = \mZ$, so that every $(\iota , \rho^s)$ is a Hilbert $\mZ$-bimodule in $\mA$.

A Hilbert $\mB$-bimodule $\mM$ in $\mA$ is said to be {\em finitely generated} if there is a finite set $\left\{ \psi_l \right\} \subset \mM$ such that for every $\psi \in \mM$ we find $\psi = \sum_l \psi_l \psi_l^* \psi$ (so that $\psi$ is a right linear combination with coefficients in $\mB$ of the $\psi_l$'s). In such a case there are natural identifications of $\mM^*$ with the dual of $\mM$, $\mM \cdot \mM^*$ with the \sC algebra $\mL(\mM)$ of right $\mB$-module endomorphisms on $\mE$, and $\mM^r := \mM \cdots \mM$ with the $r$-fold tensor power of $\mM$.  Now, the following properties hold:

\begin{itemize}

\item $P_\mM := \sum_l \psi_l \psi_l^*$ is a projection, and does not depend on the choice of the generators. We call $P_\mM$ the {\em support of} $\mM$. Of course $P_\mM \psi = \psi$ for every $\psi \in \mM$, so that $P_\mM$ corresponds to the identity map in $\mM \cdot \mM^* \simeq \mL(\mM)$.

\item The map $\sigma_\mM (a) := \sum_l \psi_l a \psi_l^*$, $a \in \mB' \cap \mA$ does not depend on the choice of the generators $\psi_l's$, and defines an endomorphism on $\mB' \cap \mA$. In a such a case we say that $\sigma_\mM$ is {\em induced by} $\mM$. Note that if $\mA$ is unital then $\sigma_\mM (1) = P_\mM$. It is clear that if $\mB$ is contained in $\mZ$ then $\sigma_\mM$ is an endomorphism of $\mA$.

\end{itemize}

\bigskip

Let now $M(\mA)$ denote the multipliers algebra of $\mA$, and $i : \mA \hra M(\mA)$ be the canonical immersion. Note that the minimal unitization $\mZ^+$ of $\mA$ appears in $M(\mA)$ as a unital \sC subalgebra of the centre $ZM(\mA)$. We say that an endomorphism $\rho$ of $\mA$ is {\em inner} if there exists a non zero, finitely generated Hilbert $\mZ^+$-bimodule $\mM \subset M(\mA)$ such that $i \circ \rho = \sigma_\mM \circ i$. Note that if $\mA$ is unital then $\mA = M(\mA)$, $\mZ = \mZ^+$, and $\mM$ is uniquely recovered as the space of intertwiners $(\iota , \rho)$. 

Our definition naturally generalizes the usual notion of inner endomorphism; in fact if $s$ is an isometry in $M(\mA)$ defining the endomorphism $\rho : \rho (a) = s a s^*$ then $(\iota , \rho) = ZM(\mA) \cdot s$. Note that we can multiply $s$ by unitaries of $ZM(\mA)$ and get the same endomorphism $\rho$. Thus the object canonically associated to $\rho$ is the bimodule $ZM(\mA) \cdot s$ rather then the isometry $s$.

\begin{ex}
Let $\mA$ be a continuous trace \sC algebra with spectrum $X$. A locally unitary automorphisms $\alpha$ in the sense of \cite{PR84} defines a line bundle over $X$. The associated module of continuous sections $\mM_\alpha$ is then a Hilbert $C_b(X)$-bimodule in the multipliers algebra $M(\mA)$.
\end{ex}

\begin{ex}
Let $\coe$ be the Pimsner algebra associated to the Hilbert $C(X)$-bimodule $\wE$ of continuous sections of a vector bundle $\mE$ over a compact space $X$ (we set left $C(X)$-action coinciding with the right one). Then $\wE$ is a Hilbert $C(X)$-bimodule in $\coe$, and induces the canonical shift endomorphism (see \cite{Vas,Vas01}).
\end{ex}

\section{Crossed Products and Vector Bundles}
\label{cp_vec_bundle}

Let $\rho$ be an endomorphism of a \sC algebra $\mA$. Then $\rho$ induces a natural structure of $C(X)$-algebra on $\mA$. In fact if $M(\mA)$ is the multipliers algebra of $\mA$ with centre $ZM(\mA)$ we define

\begin{equation}
\label{def_zro}
\zro := \left\{ f \in ZM(\mA) : \rho (fa) = f \rho(a) , a \in \mA  \right\} .
\end{equation}

\noindent Being $\zro$ a unital subalgebra of $ZM(\mA)$ it is clear that $\mA$ is a $\zro$-algebra. Furthermore by definition $\rho$ is equivariant w.r.t. the $\zro$-action. Note also that $\zro$ is an invariant of the sector 

\begin{equation}
\label{sector}
{\bf sect}(\rho) := \left\{ \sigma \in {\bf end}\mA : \sigma (a) = u \rho (a) u^* , uu^* = u^*u = 1 , u \in M(\mA) \right\}.
\end{equation}

The previous elementary remark suggests to consider $C(X)$-algebras as a natural framework for endomorphisms and related crossed product constructions.

\bigskip

Let now $X$ be a locally compact Hausdorff space, $\mA$ a $C_0(X)$-algebra, $\rho$ a $C_0(X)$-endomorphism on $\mA$. Given a vector bundle $\mE$ over $X$ with compact support (i.e. $\mE$ is trivial in a neighborhood of the point at infinity: in other terms, $\mE$ defines a vector bundle over the one point compactification $X^+$), we want to construct a crossed product of $\mA$ by $\rho$ in such a way that the endomorphism becomes induced by the module $\wE$ of continuous sections of $\mE$. As particular cases we obtain the ordinary crossed product constructions \cite{Cun82},\cite{Pas80},\cite{Sta93}, recovered by choosing trivial vector bundles eventually with rank $1$ (see \cite{Vas01} for details).

Before to proceed recall that every $C_0(X)$-algebra $\mA$ can be faithfully represented over a suitable Hilbert $C_0(X)''$-bimodule $\mM$ by means of a $C_0(X)$-monomorphism (see \cite{Kas88}). Note that $C(X^+)$ is contained in the centre of the \sC algebra $\mL(\mM)$ of right $C_0(X)''$-module endomorphisms of $\mM$.

\begin{defn}
A covariant representation of $\mA,\rho$ with rank $\mE$ is a pair $(\pi , \mE_\pi)$, where 

\begin{itemize}

\item $\pi : M(\mA) \ra \mL(\mM_\pi)$ is a (unital) non degenerate $C_0(X)$-representation of the multipliers algebra $M(\mA)$ over a Hilbert $C_0(X)''$-module $\mM_\pi$;

\item $\mE_\pi$ is a Hilbert $C(X^+)$-bimodule in $\mL(\mM_\pi)$ isomorphic to $\wE$, and such that the relation $\pi \circ \rho = \sigma_{\pi} \circ \pi$ holds; $\sigma_\pi$ is the inner endomorphism on $\mL (\mM)$ induced by $\mE_\pi$.

\end{itemize}
\end{defn}

We are going to construct our crossed product in such a way that a universality condition is satisfied w.r.t. covariant representations. As a first step we prove the existence of such representations. At this purpose let us introduce the inductive limit $\mA_\infty := \mA \stackrel{\rho}{\ra} \mA \stackrel{\rho}{\ra} \cdots$.

\begin{lem}
Let $\mE$ be a vector bundle over $X$ with compact support. There exist covariant representations $(\pi , \mE_\pi)$ of $\mA,\rho$ with rank $\mE$ if and only if $\mA_\infty \neq \left\{ 0 \right\}$.
\end{lem}

\begin{proof}
We proceed as in \cite{Sta93} and related references. Let $\mA_\infty \neq \left\{ 0 \right\}$. We consider the crossed product $\cpn$ of $\mA$ by $\rho$ with canonical maps $\mA \stackrel{\rho_0}{\ra} \mA_\infty \stackrel{i_\infty}{\hra} \cpn$ and define $i_\mA := i_\infty \circ \rho_0$. Recall that $\cpn$ is obtained by the crossed product $\mA_\infty \rtimes_{\rho_\infty} \bZ$ of $\mA_\infty$ by the automorphism $\rho_\infty$ induced by $\rho$, with canonical injection $j_\infty : \mA_\infty \hra \mA_\infty \rtimes_{\rho_\infty} \bZ$, as the corner $\cpn := p_\infty \cdot \mA_\infty \rtimes_{\rho_\infty} \bZ \cdot p_\infty$, where $p_\infty := j_\infty \circ \rho_0 (1)$ (here $1$ denotes the identity of $M(\mA)$, note that by surjectivity $\rho_0$ extends naturally on the multipliers algebra). Let now us prove that $\cpn$ is a $C_0(X)$-algebra. As $\mA_\infty \neq \left\{ 0 \right\}$ we find that $\rho_0 (a)$ is non zero for some $a \in \mA$. Now for every $f \in C_0(X)$ we find $\rho^k(fa) = f \rho^k(a) \neq 0$ thus, regarding at this identity in $\mA_\infty$, we find $\rho_\infty \circ \rho_0 (fa) = \rho_0 (f) \cdot \rho_\infty \circ \rho_0 (a)$. Let now $v \in \cpn$ be the canonical isometry such that $i_\mA \circ \rho (a) = v i_\mA(a) v^*$, with $p := vv^*$. The relation $\rho (fa) = f \rho(a)$ implies $i_\mA (f) v i_\mA (a) v^* =  v i_\mA (f) i_\mA (a) v^*$, so that (being $i_\mA$ non degenerate) $i_\mA(f)$ commutes with $i_\mA(a)$, $v$. Thus $\cpn$ is a $C_0(X)$-algebra. By \cite{Kas88} we know that there exist a faithful $C_0(X)$-module representation $\nu$ of $\cpn$ over a Hilbert $C_0(X)''$-module $\mM$. Furthermore we consider a faithful $C(X^+)$-module representation $\nu'$ of $\coe$ over a Hilbert $C_0(X)''$-module $\mM'$. Now (following \cite{Kas88}, Def. 1.6) we consider the (ungraded) tensor product $\mM \otimes_{C_0(X)''} \mM'$ and claim that the pair $\left( (\nu \circ i_\mA) \otimes 1 , \nu(v) \otimes \nu' (\wE) \right)$ is a covariant representation of $\mA$ over $\mM \otimes_{C_0(X)''} \mM'$. In fact it is obvious that $\nu(v) \otimes \nu' (\wE)$ is isomorphic to $\wE$ as a Hilbert $C(X^+)$-bimodule, and has support

$$
\sum \limits_l ( \nu(v) \otimes \psi_l) \cdot ( \nu(v)^* \otimes \psi_l^*) = (\nu (v) \cdot \nu(v^*)) \otimes 1 = \nu (p) \otimes 1 ,
$$

\noindent where $\left\{ \psi_l \right\}$ is a finite set of generators of $\nu' (\wE)$. If $a \in \mA$, $\psi \in \nu'(\wE)$ then

$$
\begin{array}{ll}

\left( \nu \circ i_\mA (\rho(a)) \otimes 1 \right)  \cdot  \left( \nu(v) \otimes \psi \right) & = 

\nu \left(  i_\mA (\rho(a)) \cdot v   \right)  \otimes \psi  = \\ & =

\nu \left(  i_\infty \circ \rho_\infty (\rho_0(a)) \cdot v   \right)  \otimes \psi   = \\ & =

\nu \left(  v \cdot i_\mA(a) \cdot v^* \cdot v   \right)  \otimes \psi   = \\ & =

\nu \left(  v \cdot i_\mA(a)  \right)  \otimes \psi   = \\ & =

\left( \nu ( v )  \cdot \nu \circ i_\mA(a) \right)  \otimes \psi   = \\ & =

\left( \nu ( v ) \otimes \psi \right)  \cdot \left( \nu \circ i_\mA(a) \otimes 1 \right) 

\end{array}
$$

\noindent and the first implication is proved. Viceversa if $\mA_\infty = \left\{ 0 \right\}$ then for every $a \in \mA$ there exists a $k \in \bN$ such that $\rho^k(a) = 0$. Thus if $(\pi , \mE_\pi)$ is a covariant representation then $0 = \pi(\rho^k(a)) = \sum \limits_{ \left| L \right| = k } \psi_L \pi (a) \psi_L^*$ which implies $\pi (a) = 0$, because $\pi \circ \rho^k$ has left inverse \footnote{We recall that a left inverse of an endomorphism $\rho$ is a positive linear map $\phi : \mA \ra \mA$ such that $\phi(1) = 1$ and $\phi(a \rho(b)) = \phi(a)b$ for every $a,b \in \mA$.} the map $a \mapsto  \frac 1{d^k} \sum \limits_L \psi_L^* \pi(a) \psi_L$, where $d$ is the rank of $\mE$. This gives a contradiction, so that the lemma is proved.

\end{proof}

\begin{prop}
Let $\mA$ be a $C_0(X)$-algebra and $\rho$ a $C_0(X)$-endomorphism on $\mA$ with $\mA_\infty \neq \left\{ 0 \right\}$. Then for every vector bundle $\mE$ on $X$ with compact support there exists up to isomorphism a unique $C_0(X)$-algebra $\cpen$ such that

\begin{itemize}

\item there is a non degenerate $C_0(X)$-homomorphism $i : \mA \ra \cpen$;

\item The module of continuous sections $\wE$ of $\mE$ is contained in $M(\cpen)$ as a finitely generated Hilbert $C(X^+)$-bimodule, and $i \circ \rho = \sigma_\mE \circ i$, where $\sigma_\mE$ is the endomorphism induced by $\wE$;

\item $\cpen$ is generated as a \sC algebra by $i(\mA)$ and $\wE$;

\item for every covariant representation $(\pi , \mE_\pi)$ with rank $\mE$ there is a non degenerate $C_0(X)$-representation $\Pi : \cpen \ra \mL(\mM_\pi)$ such that $\Pi \circ i = \pi$ and $\Pi (\wE) = \mE_\pi$.

\end{itemize}

\end{prop}

\begin{proof}
We consider a set $\Gamma$ of covariant representations such that every covariant representation of $\mA$ is equivalent to a unique direct sum of elements of $\Gamma$. Then we define the representation $\pi_\Gamma := \oplus_{\gamma \in \Gamma} \gamma$. Now, given the isomorphism $i_\gamma : \wE \ra \mE_\gamma$ associated to the covariant representation $\gamma$, we define $\mE_\Gamma := \left\{ \oplus_\gamma i_\gamma(\psi) , \psi \in \wE  \right\}$. $\mE_\Gamma$ is isomorphic to $\wE$ as a Hilbert $C(X^+)$-module and it is easily verified that $(\pi_\Gamma , \mE_\Gamma)$ is a covariant representation. Thus we introduce the \sC algebra generated by $\pi_\Gamma (\mA) \cdot \mE_\Gamma$, that by construction is a crossed product as required. Unicity up to isomorphism follows by universality.
\end{proof}

\begin{ex}
\label{ex_coe}
Let $\coe$ be the Pimsner algebra associated to the $C(X)$-bimodule of continuous sections of a vector bundle $\mE$ over a compact space $X$ (with coinciding left and right actions). Let furthermore $\coe^0 \simeq \lim_r (\mE^r,\mE^r)$ be the fixed point \sC algebra of $\coe$ w.r.t. the canonical circle action (here $\mE^r$ denotes the $r$-fold tensor power of $\mE$ and $(\mE^r,\mE^r)$ is the space of bimodule maps from $\mE^r$ into $\mE^r$). If $p \in (\mE,\mE)$ is a projection we define the shift endomorphism $\wa p (t) := t \mapsto p \otimes t , t \in (\mE^r,\mE^r)$. In particular if $p$ is a rank $1$ projection then $\coe \simeq \coe^0 \rtimes_{\wa p}^\mL \bN$, where $\mL := p \mE$ is the line bundle associated to $p$ (details in \cite{Vas01}).
\end{ex}

\section{$C(X)$-categories.}
\label{categories}

Let $\mC$ be a category with objects $\rs,\cdots$. Then $\mC$ is said to be a \sC {\em category} (see \cite{DR89} for an exhaustive introduction) if every space of arrows $(\rs)$ is a Banach space, composition of arrows is norm decreasing and there is an isometric involutive cofunctor $* : \mC \ra \mC$ satisfying the \sC identity w.r.t. composition. Functors on \sC categories preserving the above structure are called \sC {\em functors}. Note that a \sC category with a single object is just a unital \sC algebra, given by the unique space of arrows. A {\em strict tensor} \sC category (\cite{DR89}) is a \sC category with an associative bilinear \sC bifunctor (called the {\em tensor product}) having a unit (the so called identity object) and commuting with $*$. Basic examples of strict tensor \sC categories are the duals of compact groups, endowed with the tensor product of unitary representations over finite dimensional Hilbert spaces.

\bigskip

Let now $\rho$ be an endomorphism of a unital \sC algebra $\mA$. We consider the category $\wa \rho$ with objects the positive integers $r,s \in \bN$ and arrows the intertwiners $(\rho^r,\rho^s)$. Multiplication and involution of $\mA$ induce a structure of \sC category on $\wa \rho$, that comes also equipped with the associative tensor product 

\begin{equation}
\label{tens_end}
\left\{ 
\begin{array}{l}
r,s \mapsto r + s  \\ 
t,s \mapsto t \rho^r (t') = \rho^s (t') t \in (\rho^{r+r'} , \rho^{s+s'})
\end{array}
\right.
\end{equation}

\noindent where $t \in (\rho^r,\rho^s) , t' \in (\rho^{r'},\rho^{s'})$. The eventuality that $\wa \rho$ is symmetric (in the sense of \cite{DR89}) is encoded in the following property of $\rho$.

\begin{defn}[Permutation symmetry]
\label{def_perm_sim}
A unital endomorphism $\rho$ of a \sC algebra $\mA$ has permutation symmetry if there is a unitary representation $p \mapsto \eps (p)$ of the group $\bP_\infty$ of finite permutations of $\bN$ in $\mA$ such that:

\begin{equation}\label{ps1} \eps (\bS p) = \rho \circ \eps (p) \end{equation}
\begin{equation}\label{ps2} \eps := \eps (1,1) \in (\rho^2 , \rho^2) \end{equation}
\begin{equation}\label{ps3} \eps (s,1) t = \rho (t) \eps (r,1), t \in (\rho^r , \rho^s) \end{equation}

\noindent where $(r,s) \in \bP_{r+s}$ permutes the first $r$ terms with the remaining $s$ and $\bS$ is the shift $(\bS p)(1) := 1$, $(\bS p)(n) := 1 + p(n-1)$, $p \in \bP_\infty$.
\end{defn}

\noindent In the case in which $\mA$ has non trivial centre the above definition is not satisfied in interesting particular cases, as for example the shift endomorphism defined on suitable subalgebras of Pimsner algebras (see \cite{DPZ97,Vas01}). In particular, not all intertwiners in $(\rho^r,\rho^s)$ may satisfy (\ref{ps3}). We say that $\rho$ has {\em weak permutation symmetry} if just (\ref{ps1},\ref{ps2}) are satisfied and denote by $[\rho]$ the strict tensor \sC subcategory of $\rho$ having arrows

$$
[\rho^r,\rho^s] := \left\{ t \in (\rho^r,\rho^s) : \eps (s,1) t = \rho (t) \eps (r,1) \right\}.
$$

\noindent Note that $[\rho^r,\rho^r]$ certainly contains the unitaries in $\eps(\bP_r)$ if $\rho$ carries weak permutation symmetry. We denote by $\mO_{[\rho]}$ the \sC algebra generated by the $[\rho^r,\rho^s]$'s. It is easily verified that $\mO_{[\rho]}$ is a $\rho$-stable $\zro$-algebra.

\bigskip

The strict tensor \sC category $\wa \rho$ has been deeply investigated by Doplicher and Roberts. In the trivial centre case $\mZ := \mA \cap \mA' = \bC 1$, one of key results is the characterization of endomorphisms $\rho$ admitting an embedding of \sC tensor categories $\wa \rho \hra \wa \mH$, where $\wa \mH$ is the category having as object the tensor powers of a finite dimensional Hilbert space $\mH$, equipped with the usual tensor product (see \cite{DR89A}).

In the case in which $\mA$ has non trivial centre further structure is carried by $\wa \rho$. Aim of the present section is to introduce such further categorical structures.

\bigskip

Let $\mC$ be a \sC category. Recall that for every pair of objects $\rs$ of $\mC$, composition of arrows induces on the space of intertwiners $(\rs)$ a natural structure of Hilbert bimodule, acted on the left (resp. right) by the \sC algebra $(\sss)$ (resp. $(\rr)$). A family of Banach space maps

$$
F := \left\{ F_{\rs} : (\rs) \ra (\rs) \right\}_{\rs \in \mC}
$$ 

\noindent is said to be a right (resp. left) multiplier of $\mC$ if the following conditions are satisfied:

\begin{itemize}

\item  F is uniformely bounded, i.e. $\sup_{\rs} \left\| F_{\rs} \right\| < \infty$;

\item $F_{\rho,\tau}(t \circ t') = t \circ F_{\rs}(t')$ (resp. $F_{\rho,\tau}(t \circ t') = F_{\sigma,\tau}(t) \circ t'$), where $t' \in (\rs)$, $t \in (\sigma, \tau)$ .

\end{itemize}

\noindent Note that left and right multipliers can be composed in the natural way $F \circ G := \left\{ F_{\rs} \circ G_{\rs} \right\}$, and that every left (right) multiplier defines on each $(\rs)$ a right Hilbert $(\rr)$-module map (resp. a left Hilbert $(\sss)$-module map). A {\em multiplier} is a pair $F := (F^l,F^r)$, where $F^r$ (resp. $F^l$) is a right (left) multiplier, satisfying the relation

\begin{equation}
\label{rel_mult}
t \circ F^l_{\rs}(t') = F^r_{\sigma,\tau}(t) \circ t'
\end{equation}

\noindent for each $t,t'$ as above. We denote by $M(\mC)$ the set of multipliers on $\mC$ and introduce the notation $M_{\rs}(\mC) := \left\{ F_{\rs} : F \in M(\mC) \right\}$.

An obvious example is given by elements $a \simeq (a_\rho) \in \oplus_\rho (\rr)$, that define left (right) multipliers $a^l_{\rs} (t) := a_\sigma \circ t$ (resp. $a^r_{\rs} (t) := t \circ a_\rho$).

Note that being each $(\rr)$ a unital \sC algebra then $M_{\rr}(\mC) = (\rr)$. $M(\mC)$ has an obvious structure of vector space; furthermore in analogy with the multipliers of a \sC algebra we can define the *-algebraic structure $F \cdot G := (F^l \circ G^l , G^r \circ F^r)$ and $F^* := ((F^r)^*,(F^l)^*)$, where $\left\{ (F^{l,r}_{\rs})^* (t) := (F^{r,l}_{\sigma,\rho} (t^*))^* \right\}$. It is clear that $(F^l_{\rs})^* , (F^r_{\rs})^*$ are the adjoints operator of $F^l_{\rs} , F^r_{\rs}$ on the Hilbert bimodules $(\rs)$.

\begin{prop}
$M(\mC)$ has a natural structure of unital \sC algebra, and there is an isomorphism $\oplus_{\rho \in \mC} (\rr) \simeq M(\mC)$. If $Z \in M(\mC)$ satifies the relation $Z^l = Z^r$ then belongs to the center of $M(\mC)$. Furthermore \sC functors $\Phi : \mC \ra \mC'$ induce homomorphisms $\Phi_* : M(\mC) \ra M(\mC')$.

\end{prop}

\begin{proof}
The \sC algebra structure of $M(\mC)$ and the inclusion of $\oplus_\rho (\rr)$ have been proved in the above remarks. Viceversa if $F \in M(\mC)$ and $t \in (\rs)$ then $F^l_{(\rs)} (t) = F^l_{(\rs)} (1_\sigma \circ t) = F^l_{(\sss)} (1_\sigma) \circ t$. Let now $Z$ such that $Z^l = Z^r$. Then for every multiplier $F$ and $t \in (\sigma,\tau)$, $t' \in (\rs)$, we find 

$$
\begin{array}{ll}

F^r_{\sigma,\tau} \cdot Z^r_{\sigma,\tau} (t) \circ t' & = 
Z^r_{\sigma,\tau} (t) \circ F^l_{\rs} (t')  = \\ & =
Z^l_{\sigma,\tau} (t) \circ F^l_{\rs} (t')  = \\ & =
Z^l_{\sigma,\tau} (t \circ F^l_{\rs} (t'))  = \\ & =
Z^l_{\sigma,\tau} (F^r_{\rs} (t) \circ t')  = \\ & =
Z^l_{\sigma,\tau} \cdot F^r_{\sigma,\tau} (t) \circ t'  = \\ & =
Z^r_{\sigma,\tau} \cdot F^r_{\sigma,\tau} (t) \circ t' .
\end{array}
$$

\noindent Thus $F^r Z^r = Z^r F^r$. In the same way we find $F^l Z^l = Z^l F^l$, so that $Z$ belongs to the center of $M(\mC)$. 
\end{proof}

We denote by the notation $Z(\mC)$ the abelian \sC subalgebra of $M(\mC)$ of multipliers satisfying the relation $Z^l = Z^r$. In the case in which $\mC$ has a single object (i.e. $\mC$ is a unital \sC algebra) then $Z(\mC)$ is simply the centre of $\mC$.

\begin{defn}
Let $X$ be a compact Hausdorff space. A \sC category $\mC$ is said to be a $C(X)$-category if there is a unital homomorphism $C(X) \ra Z(\mC)$. $C(X)$-functors between $C(X)$-categories are \sC functors equivariant w.r.t. the $C(X)$-actions.
\end{defn}

The choice of considering compact spaces is dictated by the fact that the space of arrows of a given object in a \sC category has to be a {\em unital} \sC algebra. Now by definition if $\rs$ is a couple of objects in $\mC$ then $(\rr)$ is a $C(X)$-algebra, and $(\rs)$ is a $C(X)$-Hilbert $(\sss)$-$(\rr)$-module in the sense of \cite{Kas88,Lan01,PT00}. Thus for every $\rs$ objects in $\mC$ we can assign the "fibers" over $x \in X$ by defining 

$$
(\rs)_x := (\rs) \otimes_{(\rr)} (\rr)_x ,
$$

\noindent where $(\rr)_x$ is the quotient of $(\rr)$ by the ideal $(\rr) \cdot \mC_x(X)$ (here $\mC_x(X)$ denotes the ideal of continuous functions vanishing on $x$). We have the following categorical analogue of the Nilsen picture of $C(X)$-algebras:

\begin{prop}
Let $X$ be a compact Hausdorff space. If $\mC$ is a $C(X)$-category then there exists a family $\left\{ x_* : \mC \ra \mC_x \right\}_{x \in X}$ of \sC functors such that for every arrow $t \in (\rs)$ the following properties hold, with $t_x := x_*(t)$:

\begin{itemize}
\item the norm function $n(t) : x \mapsto \left\| t_x \right\|$ is upper semicontinuous;
\item $\left\| t \right\| = \sup_x \left\| t_x \right\|$.
\end{itemize}
\end{prop}

\begin{proof}
Let $x \in X$. We consider the category $\mC_x$ having as objects the pairs $\rho_x := \left\{ \rho , x \right\}$, where $\rho$ is an object of $\mC$, and arrows the Banach spaces $(\rho_x , \sigma_x) := (\rs)_x$. If $t \in (\rs)$ we define $t_x \in (\rho_x , \sigma_x)$ to be the image of $t$ by the canonical projection $(\rs) \ra (\rs)_x$. Composition in $\mC_x$ is naturally given by $t_x \circ t'_x := (t \circ t')_x$, and is well defined because of the fact that composition in $\mC$ is a $C(X)$-bimodule map. Now $\left\| t \right\|^2 = \left\| t^* \circ t \right\|$, where $t^* \circ t$ belongs to the $C(X)$-algebra $(\rr)$. So that by \cite{Nil96} the norm function $n(t)$ is upper semicontinuous, with $\sup$ exactly the norm of $t$. Thus we have defined for $x \in X$ the \sC functor $x_*$, by assigning $\rho \mapsto \rho_x$ and $t \mapsto t_x$ for $t \in (\rs)$.
\end{proof}

We summarize the properties exposed in the previous proposition by saying that $\mC$ is an {\em upper semicontinuous field of C*-categories over} $X$, and call {\em fiber functors} the $x_*$'s. The previous result implies that the usual operation of restriction of a $C(X)$-algebra on a compact subset $Y \subset X$ makes sense for $C(X)$-categories. This fact allows us to give the notion of local isomorphism:

\begin{defn}
$C(X)$-categories $\mC$, $\mC'$ are said to be locally isomorphic if for every $x \in X$ there exist a compact neighborhood $U$ of $x$ such that the restrictions $\mC_U$, ${\mC'}_U$ are isomorphic as $C(U)$-categories.
\end{defn}

Given a \sC category $\mC$, we can consider the {\em constant field} $C( X ,\mC )$ having the same objects of $\mC$ and arrows the spaces of continuous functions $C (X , (\rs) )$, where $\rs$ are objects in $\mC$.

\begin{defn}
A $C(X)$-category $\mC$ is said to be a continuous field of \sC categories if for each $t \in (\rs)$ arrow in $\mC$ the norm function $n(t)$ is continuous. $\mC$ is said to be (locally) trivial if it is (locally) isomorphic as a $C(X)$-category to a constant field of \sC categories.
\end{defn}

Thus a locally trivial $C(X)$-category is a continuous field of \sC categories. In particular, if $\rho$ is an object of a continuous field of \sC categories then $(\rr)$ is a continuous field of \sC algebras.

\bigskip

A $C(X)$-category endowed with a $C(X)$-tensor product (which is denoted by $\times$) is said to be a {\em strict tensor} $C(X)$-{\em category}. We now prove that every strict tensor \sC category has a natural structure of strict tensor $C(X)$-category. We start by introducing the following abelian \sC algebra, given a strict tensor \sC category $\mC$ with identity object $\iota$:

\begin{equation}
[\ii] := \left\{ f \in (\ii) : 1_\rho \times f = f \times 1_\rho , \rho \in {\bf obj}(\mC)  \right\}
\end{equation}

\noindent In the case $[\ii] = \bC 1$ we say that $\mC$ is an {\em ergodic} category.

\begin{prop}
Let $\mC$ be a strict tensor \sC category with identity object $\iota$. Then $\mC$ is a strict tensor $[\ii]$-category.
\end{prop}

\begin{proof}
First we prove that $\mC$ is a $[\ii]$-category. We define the map $[\ii] \ra Z(\mC)$ by posing $f \mapsto (f^l , f^r) : f^l_{\rs}(t) := f \times t , f^r_{\rs}(t) := t \times f$, where $t \in (\rs)$. As $f \in [\ii]$ we find $f^r = f^l$, and $(f^l , f^r)$ belongs to $Z(\mC)$. Furthermore the map $f \mapsto (f^l , f^r)$ is obviously unital and thus non degenerate, so that $\mC$ is a $[\ii]$-category. In the sequel we will write simply $ft := f^l_{\rs} (t)$. We now have to prove that the tensor product is equivariant w.r.t the $[\ii]$-action, but this fact is obvious because of the relation $f \times t \times s' = ft \times t' = t \times f t'$, $f \in [\ii], t \in (\rs) , t' \in (\rho',\sigma')$, and analogues.
\end{proof}

\begin{cor}
Every strict tensor \sC category defines an upper semicontinuous field of tensor \sC categories over the spectrum of $[\ii]$.
\end{cor}

We now prove a simple result on the structure of the DR-algebra associated to an object $\rho$ of $\mC$. As customary for simplicity we assume that the operation of tensoring on the right with the identity arrow $1_\rho \in (\rr)$ is an injective map.

\begin{lem}
Let $\mC$ be a strict tensor $C(X)$-category (resp. a continuous field of strict tensor \sC categories). Then for each $k \in \bZ$ the Banach space $\orok := ...\stackrel {\times 1_\rho} \longrightarrow (\rho^r,\rho^{r+k}) \stackrel {\times 1_\rho} \longrightarrow ( \rho^{r+1},\rho^{r+k+1}) \stackrel {\times 1_\rho } \longrightarrow ...$ has the structure of an upper semicontinuous (resp. a continuous) field of Banach spaces $(\orok , \orok \ra \mO_{\rho_x}^k )$).
\end{lem}

\begin{proof}
It suffice to observe that the embeddings 

$$...\stackrel {\times 1_\rho} \longrightarrow
(\rho^r,\rho^{r+k}) \stackrel {\times 1_\rho} \longrightarrow
( \rho^{r+1},\rho^{r+k+1}) \stackrel {\times 1_\rho } \longrightarrow ...
$$ 

\noindent given by tensoring on the right by $1_\rho$ are $C(X)$-module maps. Thus there is a dense subset of $\orok$ whose elements have upper semicontinuous (resp. continuous) norm function, and $\orok$ is an upper semicontinuous (resp. continuous) field of Banach spaces.
\end{proof}

In particular, note that $\mO_\rho^0$ is a $C(X)$-algebra (resp. a continuous field of \sC algebras). We denote by $\wa \rho$ the strict tensor $C(X)$-category generated by the tensor powers of $\rho$.

\begin{prop}
Let $\rho$ be an object of a strict tensor $C(X)$-category. Then $\oro$ is a $C(X)$-algebra, and is a (locally trivial) continuous field $(\oro , \oro \ra \mO_{\rho_x})$ over $X$ if $\wa \rho$ is a (locally trivial) continuous field itself.
\end{prop}

\begin{proof}
As there is a unital monomorphism of $C(X)$ into the centre of $\oro$ we obtain that $\oro$ is a $C(X)$-algebra. Thus $\oro$ defines an upper semicontinuous field over $X$. We have now to prove that norm function is continuous for elements of $\oro$ if $\wa \rho$ is a continuous field of \sC categories. At this purpose note that by previous lemma there is a dense set in $\oro$ of operators having continuous norm function (belonging to the $\oro^k$'s), so that $\oro$ is a continuous field. Local triviality follows by funtoriality of $\oro$: in fact local charts $F_U : \left. \mC \right|_U \ra C(U,\mC_x)$ define $C(U)$-algebra homomorphisms $\left. \oro \right|_U \ra C(U) \otimes \mO_{\rho_x}$.
\end{proof}

\begin{ex}
The category having as objects endomorphisms belonging to a sector ${\bf sect}(\rho)$ (defined in (\ref{sector})), with arrows the intertwiners. In this case $C(X) = \zro$. Furthermore ${\bf sect}(\rho)$ is a strict tensor $\zro$-category with tensor product $\rho \times \sigma := \rho \circ \sigma$, $t \times t' := t \rho (t')$, $t \in (\rho,\rho')$, $t' \in (\sigma,\sigma')$, and $(\ii) = \mZ$, $[\ii] = \zro$.
\end{ex}

\section{Cross products by dual actions.}
\label{dual_dyn_sys}

Let $C(X)$ be an abelian \sC algebra. In the sequel, a Hilbert $C(X)$-bimodule with left action coinciding with the right one will be denoted by the notation $\mE$, emphasizing the fact that $\mE$ is recovered up to isomorphism by the corresponding vector bundle over $X$. It is then proved that $\coe$ carries a natural structure of continuous field of Cuntz algebras over $X$ (see \cite{Vas}); furthermore the canonical endomorphism $\sigma_\mE$ induced by the shift $t \mapsto 1 \otimes t$, $t \in (\mE^r,\mE^s)$, is well defined on $\coe$. If $G$ is a compact group acting on $\mE$ by unitary endomorphisms then a natural action on $\coe$ is induced, as for Cuntz algebras, by extending the map $\psi \mapsto g \psi$, $\psi \in \mE$. We denote by $\coe^G$ the fixed point algebra, that is generated by the intertwiners spaces $(\mE^r,\mE^s)_G := \left\{ t \in (\mE^r,\mE^s) : g^{\otimes^s} \cdot t = t \cdot g^{\otimes^r} , g \in G \right\}$ of tensor powers of the defining representation of $G$. It is easily checked that $\coe^G$ is $\sigma_\mE$-stable. We also denote by ${\wa G}_\mE$ the strict tensor $C(X)$-category having objects the tensor powers $\mE^r$, with arrows $(\mE^r,\mE^s)_G$, endowed with the usual tensor product of Hilbert bimodule maps.

In the case in which left $C(X)$-action over an Hilbert bimodule $\mM$ does not coincide with the right one some structural properties change. In particular the canonical endomorphism $\sigma_\mM$ is just defined on the \sC subalgebra $\mB_\mM \subset \mO_\mM$ generated by $C(X)$-bimodule maps in $(\mM^r,\mM^s)$ commuting with left $C(X)$-action. If a compact group $G$ acts on $\mM$ in such a way that the fixed point algebra $\mO_\mM^G$ is contained in $\mB_\mM$, then $\mO_\mM^G$ is $\sigma_\mM$-stable. In such a case ${\wa G}_\mM$ defined as above is a strict tensor $C(X)$-category.

\bigskip

Let now $(\mB,G)$ be a \sC dynamical system. We suppose that $G$ is compact and $\mB$ is unital. Furthermore we introduce the notations $\mA := \mB^G$, $\mZ := \mA \cap \mA'$. The following notion is related to the one of Hilbert \sC system (see \cite{BL97} an related references).

\begin{defn}
With the above notations, $(\mB,G)$ is said to be a {\em cross product by a dual action} if there exists a finitely generated $G$-stable Hilbert $(\mB \cap \mB')^G$-bimodule $\mE \subset \mB$ with support $1$ such that $\mB$ is generated by $\mA , \mE$ as a \sC algebra. We then use the notation $(\mB,G,\mA,\mE)$.
\end{defn}

Note that in such a case there is an inclusion $\coe \subset \mB$, where $\coe$ is the Pimsner algebra associated to $\mE$. Furthermore $\mE$ induces an inner endomorphism $\sigma_\mE$ on $\mB$. Since $\sigma_\mE$ does not depend on the choice of the generators $\mA$ is $\sigma_\mE$-stable (in fact if $\left\{ \psi_l \right\}$ is a set of generators of $\mE$ the same is true for $\left\{ g \psi_l \right\}$ for $g \in G$, and $g \circ \sigma_\mE (b) = \sigma_\mE \circ g (b)$). We denote by $\rho \in {\bf end}\mA$ the restriction of $\sigma_\mE$ on $\mA$.

\begin{thm}
Let $(\mB,G,\mA,\mE)$ be a cross product by a dual action. Then the following properties hold:

\begin{itemize}

\item $\rho$ has weak permutation symmetry;

\item $\zro = (\mB \cap \mB')^G = (\mB \cap \mB') \cap \mZ$;

\item there is a functor of strict tensor $\zro$-categories ${\wa G}_\mE \ra \wa \rho$;

\item If $\mA' \cap \mB = \mZ$ then $\mM := \left\{ \psi \in \mB : \psi a = \rho (a) \psi , a \in \mA  \right\}$ is a $G$-stable Hilbert $\mZ$-bimodule in $\mB$ and $\mM = \mE \cdot \mZ$. Furthermore $\mO^G_\mM \subseteq \mB_\mM$.

\item If $\mA' \cap \mB = \mZ$ then there are isomorphisms of strict tensor $\zro$-categories ${\wa G}_\mE \simeq [\rho]$ and ${\wa G}_\mM \simeq \wa \rho$. Furthermore the following diagram of morphisms of \sC dynamical systems commutes:

\begin{center}
$\xymatrix{
            (\mO_{[\rho]} , \rho)
		    \ar@{^{(}->}[r]
		 &  (\oro , \rho)
		 \\ (\coe^G , \sigma_\mE)
		    \ar[u]^{\simeq}
		    \ar@{^{(}->}[r]
		 &  (\mO^G_\mM , \sigma_\mM)
		    \ar[u]_{\simeq}
}$ \\
\end{center}

\end{itemize}

\end{thm}

\begin{proof}
\
\begin{itemize}

\item It suffice to pick a (finite) set of generators $\left\{ \psi_l \right\}$ of $\mE$ and define the map 

$$
\eps (p) := \sum_L \psi_{l_{p(1)}} \cdots \psi_{l_{p(r)}} \psi^*_{l_r} \cdots \psi^*_{l_1} ,
$$ 

\noindent where $p \in \bP_r$. As left $(\mB \cap \mB')^G$-action on $\mE$ coincides with the right one (in fact it is simply given by left and right multiplication) the above defined map does not depend on the choice of the generators. For the same reason $\eps (p)$ is $G$-invariant and belongs to $\mA$. We leave to the reader the easy verifications of the required properties of weak permutation symmetry.

\item It is clear that $(\mB \cap \mB')^G = (\mB \cap \mB') \cap \mZ$. Furthermore $\zro$ as defined in (\ref{def_zro}) is contained in $ \mZ \subset \mA = \mB^G$, and elements of $\zro$ commutes with elements of $\mE$ (as $\zro$ is pointwise $\sigma_\mE$-invariant). Thus $\zro \subseteq (\mB \cap \mB') \cap \mZ$. Viceversa, it is clear that $(\mB \cap \mB') \cap \mZ$ is contained in $\zro$.

\item We denote by $i : \coe \hra \mB$ the natural inclusion. The $G$-action on $\mB$ restricts to the natural $G$-action on $\coe$, regarding at $G$ as a subgroup of $\mUE$. Thus we just have to prove that $i (\mE^r,\mE^s)_G \subseteq (\rho^r,\rho^s)$. Now if $t \in (\mE^r,\mE^s)$ then $t = \sum_{LM} \psi_L \psi_M^* t_{LM}$, where $\left\{ \psi_L \right\}$ generate $\mE^s$, $\left\{ \psi_M \right\}$ generate $\mE^r$, $t_{LM} \in \zro$. Furthermore by definition of the $G$-action $i(t)$ belongs to $\mA$. Thus, as $\sigma_\mE^r(b) i(\psi_L) = i(\psi_L) b$ for $b \in \mB$, we find $\rho^s(a) i(t) = \sigma_\mE^s (a) i(t) = i(t) \sigma_\mE^r (a) = i(t) \rho^r(a)$ for $a \in \mA$. So that the required functor is defined as $\mE^r \mapsto r$ on objects and by the inclusion map $i : (\mE^r,\mE^s)_G \hra (\rho^r,\rho^s)$ on arrows. As the tensor product $t \times t' = t \times 1 \cdot 1 \otimes t'$ of Hilbert bimodule map is translated as $t \sigma_\mE (t')$ in the Pimsner algebra, our functor mantains the tensor structure (note in fact that $\left. \sigma_\mE \right|_\mA = \rho$).

\item It is clear that $\mE \cdot \mZ \subseteq \mM$. Viceversa if $\psi \in \mM$ then $\psi^* \psi' \in \mA' \cap \mB = \mZ$ for every $\psi' \in \mE$, and equality is proved. Note that in general elements of $\mZ$ do not commute with elements of $\mE$. Now it is clear that $\mM$ is $\mG$-stable, and $\mO^G_\mM \subseteq \mA \subseteq \mZ' \cap \mB$. Thus $\mO^G_\mM \subseteq \mB_\mM$.

\item Spaces of arrows of ${\wa G}_\mM$ appear in $\mB$ as $(\mM^r,\mM^s)_G \simeq (\mM^s \cdot {\mM^r}^*) \cap \mA$. We denote by $j : \com \hra \mB$ the natural inclusion. Thus if $t \in (\mM^r,\mM^s)_G$ then $j(t) = \sum_{LM} \psi_L t_{LM} \psi_M^* $, where $\psi_L \in \mE^s$, $\psi_M \in \mE^r$, $t_{LM} \in \mZ$. Now, by the same argument for $i(\mE^r,\mE^s)$, we find $j (\mM^r,\mM^s)_G \subseteq (\rho^r,\rho^s)$. Viceversa, if $t' \in (\rho^r,\rho^s)$ then $t'_{LM} := \psi_L^* t' \psi_M \in \mA' \cap \mB = \mZ$ and $t' \in (\mM^s \cdot {\mM^r}^*)$. Thus $t' \in (\mM^s \cdot {\mM^r}^*) \cap \mA = j(\mM^r,\mM^s)_G$. We have proved in this way also that $(\com^G , \sigma_\mM) \simeq (\oro , \rho)$. We now prove that $j$ restricts to an isomorphism $\coe^G \ra \soro$. It is clear that $j (\mE^r , \mE^s)_G \subseteq \left[ \rho^r , \rho^s \right]$; thus we have just to verify that $\left[ \rho^r , \rho^s \right] \subseteq j (\mE^r , \mE^s)_G$ for each $r,s \in \bN$. At this purpose it suffice to observe that if $t \in \left[ \rho^r , \rho^s \right]$ then with the above notations we find $\rho (t_{LM}) = \sigma_\mE (\psi_L^*) \rho (t) \sigma_\mE (\psi_M) = \sigma_\mE (\psi_L^*) \eps (s,1) t \eps (1,r) \sigma_\mE (\psi_M) = t_{LM}$ (in fact $\eps (s,1) \psi_L = \sigma_\mE (\psi_L)$ ). Thus $t_{LM} \in \zro$, $t \in j(\coe)$ and the inclusion is proved.

\end{itemize}
\end{proof}

\begin{ex}
We refere to Ex.\ref{ex_coe}. Let $G \subseteq \mUE$. Then $G$ acts by automorphisms on $\coe$, and the system $(\coe,G,\coe^G,\mE)$ satisfies the hypothesis of previous theorem. The endomorphism induced by $\mE$ on $\coe^G$ is the shift $\wa 1$. Furthermore ${\coe^G}' \cap \coe = C(X)$ (see \cite{Vas01} for details). In particular for $G = \bT$ we get the usual circle action, with fixed point algebra the zero grade $\coe^0$; the corresponding crossed product is $(\coe,\bT,\coe^0,\mE)$.
\end{ex}

Crucial step for duality is the construction of a canonical cross product by a dual action starting from the pair $(\mA,\rho)$. At this purpose a necessary property for $\rho$ is permutation symmetry, which has to be stated in a generalized form in the non trivial centre case. Construction of the system $(\mB,G,\mA,\mE)$ has to be made by considering a suitable quotient of the crossed product introduced in \secref{cp_vec_bundle}, in such a way that the relative commutant of $\mA$ in $\mB$ reduces to the centre of $\mA$. About these aspects we refere the reader to \cite{Vas01}.

\bigskip

{\bf Acknowledgments.} The author would like to thank Professor S. Doplicher, for valuable discussions and comments.

\end{document}